\documentclass[12pt,psfig]{article}
\usepackage{graphicx,epsfig}
\usepackage{amsmath}
\usepackage{mathrsfs}
\usepackage{amssymb}
\def\bbb{\begin{eqnarray*}}

\def\eee{\end{eqnarray*}}

\begin{document}

\baselineskip=18pt
\begin{center}

\vspace{-0.6in} {\large \bf Lyapunov exponents, sensitivity, and stability\\
for non-autonomous discrete systems }
\\ [0.2in]

Hua Shao, Yuming Shi$^{*}$, Hao Zhu

\vspace{0.15in} Department of Mathematics, Shandong University \\
 Jinan, Shandong 250100, P.~R. China\\

\footnote{$^{*}$ The corresponding author.}
\footnote{ Email addresses: huashaosdu@163.com (H. Shao), ymshi@sdu.edu.cn (Y. Shi), haozhu@mail.sdu.edu.cn \\(H. Zhu).} \

\end{center}

{\bf Abstract.}  This paper is concerned with relationships of Lyapunov
exponents with sensitivity and stability for non-autonomous discrete systems.
Some new concepts are introduced for non-autonomous discrete systems,
including Lyapunov exponents, strong sensitivity at a point and in a set,
Lyapunov stability, and exponential asymptotical stability.
It is shown that the positive Lyapunov exponent at a point implies
strong sensitivity for a class of non-autonomous discrete systems.
Furthermore, the uniformly positive Lyapunov exponents in a totally
invariant set imply strong sensitivity in this set
under certain conditions. It is also shown that the negative Lyapunov
exponent at a point implies exponential asymptotical
stability for a class of non-autonomous discrete systems.
The related existing results for autonomous discrete systems
are generalized to non-autonomous discrete systems and their conditions are weaken.
One example is provided for illustration.\medskip

\medskip

{2000 {\bf \it Mathematics Subject Classification}}: 34D08; 37B55; 49K40.

\bigskip

\noindent{\bf 1. Introduction}\medskip

Chaos is a universal dynamical behavior of nonlinear systems, and becoming more and more popular in the
research of nonlinear science. Sensitivity, which is popularized by the meteorologist
Lorenz through the so-called ``butterfly effect", is widely understood as the central element of chaos. Therefore,
the study on sensitivity has attracted a lot of attention from many scholars [1, 4, 10--12, 14, 30]. However, many
results give some sufficient conditions for sensitivity, which are qualitative and thus not easy to be verified, whereas
Lyapunov exponents can characterize sensitivity in a quantitative perspective at some extent.

The concept of Lyapunov exponent was dated back to Lyapunov, when he studied the stability of solutions of ordinary
differential equations, and he proved that the solution of a regular system is stable if all the Lyapunov exponents are negative [22].
In 1946, Chetaev showed that the solution of a regular system is unstable if at least one Lyapunov exponent is positive [7].
Their works were followed by some other scholars [3, 5, 6, 8, 9, 13, 15, 16, 18--21, 23, 24]. Recently, several results about
Lyapunov instability or stability of solutions of discrete dynamical systems were also obtained. For instance, in 2004, Abraham et al.
showed that the positive Lyapunov exponents almost everywhere imply a kind of sensitivity for measure-preserving maps [2].

It was always taken for granted that positive Lyapunov exponents imply sensitivity and negative Lyapunov exponents imply
stability. However, in 2001, Demir et al. proved that this conclusion may not hold for general interval maps by two examples [9], one of which shows that a system is not sensitive at a point with positive Lyapunov exponent, and the other shows that a system is sensitive
at a point with negative Lyapunov exponent. So it is very interesting to investigate relationships of Lyapunov exponents
with sensitivity and stability. In 2010, Ko{\c c}ak and Palmer showed that the positive strong Lyapunov exponent at a point
implies sensitivity, and the negative Lyapunov exponent at a point implies Lyapunov stability under contain conditions
for differentiable interval maps [16]. Since many complex systems
occurring in the real-world problems such as physical, biological, and economical  problems
are necessarily described by
non-autonomous discrete systems,
many scientists and mathematicians focused on complexity of
non-autonomous discrete systems recently [5, 14, 17, 25--30].
Motivated by all the above works, we shall try to investigate relationships
of Lyapunov exponents with sensitivity and stability for non-autonomous discrete systems
in the present paper.

The rest of the paper is organized as follows.
Section 2 presents some basic concepts and useful lemmas.
In Sections 3 and 4,  some relationship between
the positive Lyapunov exponent and strong sensitivity,
and some relationship between negative Lyapunov exponent and exponential asymptotical stability
for non-autonomous discrete systems are investigated, respectively.
Finally, the non-autonomous logistic system is discussed in Section 5 as an illustrative example.
\bigskip

\noindent{\bf 2. Preliminaries}\medskip

In this section, some basic concepts are presented, including Lyapunov exponents, strong sensitivity, Lyapunov stability, and
exponential asymptotical stability for non-autonomous discrete systems. In addition, some useful lemmas are also presented.\medskip

We shall consider the following non-autonomous discrete system in the present paper:
\vspace{-0.2cm}$$x_{n+1}=f_n(x_n),\; n\geq0,                                                                                   \vspace{-0.2cm}\eqno(2.1)$$
where $f_n: X\to X$ is a map for each $n\geq0$, and $X$ is metric space with metric $d$.\medskip

For any $x_0\in X$, $\{x_n\}_{n=0}^{\infty}$ is called the (positive) orbit of system (2.1) starting from $x_0$ and
$x_n=f_0^n(x_0),\;n\geq0$, where $f_0^n:=f_{n-1}\circ\cdots\circ f_{0}$. Further, by $B_{\epsilon}(x_0)$ denote the open ball
of radius $\epsilon$, centered at $x_0$.\medskip

\noindent{\bf Definition 2.1.} System (2.1) is said to be sensitive at a point $x_0\in X$ if there exists a constant $\delta>0$ such that
for any neighborhood $U$ of $x_0$, there exist $y_0\in U$ and a positive integer $N$ such that $d(f_{0}^{N}(y_0),f_{0}^{N}(x_0))>\delta$,
while $\delta$ is called a sensitivity constant of system (2.1) at $x_0$. Furthermore, system (2.1) is said to be sensitive in
a nonempty set $S\subset X$ if there exists a constant $\delta>0$ such that it is sensitive at every point in $S$ with
sensitivity constant $\delta$.\medskip

\noindent{\bf Remark 2.1.} The concept of sensitivity in a nonempty set for system (2.1) is from Definition 2.3 in [25].\medskip

For any positive integer $k$, we consider another system:
\vspace{-0.2cm}$$x_{n+1}=f_n(x_n),\; n\geq k.                                                                                 \vspace{-0.1cm}\eqno(2.2)$$

\noindent{\bf Proposition 2.1.} Let $k$ be any positive integer and $f_i$ be continuous in $X$
for each $0\leq i\leq k-1$. If system (2.1) is sensitive at $x_0$ with sensitivity constant $\delta$, then
system (2.2) is also sensitive at $x_k$ with the same sensitivity constant $\delta$.\medskip

\noindent{\bf Proof.} Fix any $0<\epsilon<\delta$. By the continuity of $f_i,\;0\leq i\leq k-1$,
there exists $\epsilon'>0$ such that
\vspace{-0.2cm}$$f_{0}^{i}(B_{\epsilon'}(x_0))\subset B_{\epsilon}(x_i),\;1\leq i\leq k.                                     \vspace{-0.1cm}\eqno(2.3)$$
Since system (2.1) is sensitive at $x_0$ with sensitivity constant $\delta$, there exist $y_0\in B_{\epsilon'}(x_0)$
and a positive integer $N$ such that $d(f_{0}^{N}(y_0),f_{0}^{N}(x_0))>\delta$. It follows from (2.3) that
$N>k$. So, $d(f_{k}^{N-k}(f_{0}^{k}(y_0)),f_{k}^{N-k}(f_{0}^{k}(x_0)))=d(f_{0}^{N}(y_0),f_{0}^{N}(x_0))>\delta$.
Hence, system (2.2) is sensitive at $x_k$ with sensitivity constant $\delta$. This completes the proof.\medskip

\noindent{\bf Definition 2.2.} Let $x_0$ not be an isolated point in $X$.
System (2.1) is said to be strongly sensitive at a point $x_0$
if there exist $\delta>0$ and a neighborhood $U$ of $x_0$ such that for any given $y_0\in U$ with $y_0\neq x_0$,
$d(f_{0}^{N}(y_0),f_{0}^{N}(x_0))>\delta$ for some positive integer $N$, while $\delta$ is called a strong
sensitivity constant of system (2.1) at $x_0$. Furthermore, system (2.1) is said to be strongly sensitive in a nonempty set
$S\subset X$ without isolated points if there exists a constant $\delta>0$ such that it is strongly sensitive at every point in $S$ with
sensitivity constant $\delta$.\medskip

\noindent{\bf Definition 2.3.} System (2.1) is said to be Lyapunov stable at $x_0\in X$ if for any $\eta>0$,
there exists $\delta>0$ such that $d(f_{0}^{n}(y_0),f_{0}^{n}(x_0))\leq \eta$ for each $n\geq0$ and any $y_0\in X$
with $d(y_0,x_0)<\delta$.\medskip

\noindent{\bf Definition 2.4.} System (2.1) is said to be exponentially asymptotically stable at $x_0\in X$
if there exist positive constants $\delta$, $C$, and $\lambda$ such that $d(f_{0}^{n}(y_0),f_{0}^{n}(x_0))\leq C\exp[-\lambda n]$
for each $n\geq0$ and any $y_0\in X$ with $d(y_0,x_0)<\delta$.\medskip

Note that exponential asymptotical stability implies Lyapunov stability at $x_0$ for system (2.1)
under the assumption that $f_k$ is continuous in $X$ for each $k\geq0$; and Lyapunov stability is exactly the inverse of sensitivity.\medskip

\noindent{\bf Definition 2.5} [25, Definition 3.2]. Assume that $\{D_n\}_{n=0}^{\infty}$ and $\{E_n\}_{n=0}^{\infty}$
are two sequences of sets in $X$, and $h_n: D_n\to E_n$ is a uniformly continuous map for each $n\geq0$.
The sequence of maps $\{h_n\}_{n=0}^{\infty}$ is said to be equi-continuous in $\{D_n\}_{n=0}^{\infty}$ if
for any $\epsilon>0$, there exists $\delta>0$ such that $d(h_n(x), h_n(y))<\epsilon$ for all
$n\geq0$ and for all $x, y\in D_n$ with $d(x, y)<\delta$.\medskip

\noindent{\bf Lemma 2.1.} Let $g: E\to E'$ and $f_n: D_n\to D_n',\;n\geq0$, be maps, where $E$, $E'$, $D_n$, and $D_n'$
are nonempty sets of $X$. Assume that $g$ is uniformly continuous in $E$ and $\{f_n\}_{n=0}^{\infty}$ is equi-continuous
in $\{D_n\}_{n=0}^{\infty}$. If $\cup_{n=0}^{\infty}f_n(D_n)\subset E$, then $\{g\circ f_n\}_{n=0}^{\infty}$
is equi-continuous in $\{D_n\}_{n=0}^{\infty}$.\medskip

\noindent{\bf Proof.} Since $g$ is uniformly continuous in $E$,
for any $\epsilon>0$, there exists $\delta>0$ such that $d(g(z_{1}),g(z_{2}))<\epsilon$
for any $z_{1}, z_{2}\in E$ with $d(z_{1},z_{2})<\delta$. By the equi-continuity
of $\{f_n\}_{n=0}^{\infty}$ in $\{D_n\}_{n=0}^{\infty}$, there exists $\delta_{1}>0$ such that for each $n\geq0$ and
any $x, y\in D_n$ with $d(x,y)<\delta_{1}$, one has that $d(f_n(x),f_n(y))<\delta$, and then $d(g(f_n(x)), g(f_n(y)))<\epsilon$
since $\cup_{n=0}^{\infty}f_n(D_n)\subset E$. Hence $\{g\circ f_n\}_{n=0}^{\infty}$ is equi-continuous in $\{D_n\}_{n=0}^{\infty}$.
This completes the proof.\medskip

The above concepts are defined in general metric spaces. Next, we shall consider system (2.1) in
a non-degenerate closed interval.\medskip

In [5], Balibrea et al. extended the concept of Lyapunov exponent for a single interval map $f$ to a sequence of interval maps
$\{f_n\}_{n=0}^{\infty}$ by the formula:
\vspace{-0.2cm}$$\lambda(x_0)=\lim_{n\to\infty}\frac{1}{n}\ln|(f_{0}^{n})'(x_0)|=\lim_{n\to\infty}\frac{1}{n}\sum_{j=0}^{n-1}\ln|f_{j}'(x_j)|.\vspace{-0.3cm}$$
However, the limit does not always exist in general (see Example 5.1). Inspired by the idea given in [24]
for autonomous discrete systems, we introduce the following definition:\medskip

\noindent{\bf Definition 2.6.} Let $I$ be a non-degenerate closed interval and $f_n: I\to I$ be a $C^{1}$ map for each $n\geq0$.
The Lyapunov exponent of system (2.1) at a point $x_0\in I$ is defined by
\vspace{-0.2cm}$$\lambda(x_0):=\limsup_{n\to\infty}\frac{1}{n}\ln|(f_{0}^{n})'(x_0)|=\limsup_{n\to\infty}\frac{1}{n}\sum_{k=0}^{n-1}\ln|f_{k}'(x_k)|,
\vspace{-0.2cm}\eqno(2.4)$$
where $\{x_k\}_{k=0}^{\infty}$ is the orbit of system (2.1) starting from $x_0$.\medskip

\noindent{\bf Proposition 2.2.} The Lyapunov exponent of system (2.1) at $x_0$ equals that of system (2.2)
at $x_k$ for each $k\geq1$, where $\{x_k\}_{k=0}^{\infty}$ is the orbit of system (2.1) starting from $x_0$.\medskip

\noindent{\bf Proof.} Since the proof is trivial, its details are omitted.\bigskip

\noindent{\bf 3. Positive Lyapunov exponent implies strong sensitivity}\medskip

In this section, we shall show that the positive Lyapunov exponent implies strong sensitivity
under certain conditions.\medskip

\noindent{\bf Theorem  3.1.} {\it Let $I$ be a non-degenerate interval, $f_n: I\to I$ be a $C^{1}$ map for each $n\geq0$,
and $x_0\in I$. Assume that $\{f_{n}'\}_{n=0}^{\infty}$ is equi-continuous in $I$, and $M:=\inf\{|f_n'(x_k)|: n,k\geq0\}>0$.
If $\lambda(x_0)>0$, then system {\rm(2.1)} is strongly sensitive at $x_0$.}\medskip

\noindent{\bf Proof.} Since $\{f_{n}'\}_{n=0}^{\infty}$ is equi-continuous in $I$, $\{|f'_{n}|\}_{n=0}^{\infty}$ is equi-continuous in $I$,
and then there exists $\delta_{0}>0$ such that for each $n\geq0$ and any $x,y\in I$ with $|x-y|<\delta_{0}$,
\vspace{-0.2cm}
$$\big||f'_{n}(x)|-|f'_{n}(y)|\big|<M/2.                                  \eqno(3.1)\vspace{-0.2cm}$$
Set
$J:=\cup_{k=0}^{\infty}(x_k-\delta_{0}, x_k+\delta_{0})\cap I$.
For any $x\in J$, there exists some $k\geq0$ such that $x\in (x_k-\delta_{0}, x_k+\delta_{0})\cap I$.
By (3.1) one gets that $|f_{n}'(x)|>|f'_{n}(x_{k})|-M/2\geq M/2$ for each $n\geq0$,
which implies that $\cup_{n=0}^{\infty}\{|f_n'(x)|: x\in J\}\subset[M/2,+\infty)$.
This, together with the fact that $\{|f_{n}'|\}_{n=0}^{\infty}$ is equi-continuous in $J$ and $\ln x$ is
uniformly continuous in $[M/2,+\infty)$, yields that $\{\ln|f_{n}'|\}_{n=0}^{\infty}$ is equi-continuous in $J$ by Lemma 2.1.
So, there exists $0<\delta_{1}<\delta_{0}$ such that
for any $x,y\in J$ with $|x-y|<\delta_{1}$,
\vspace{-0.2cm}$$\big|\ln|f_{n}'(x)|-\ln|f_{n}'(y)|\big|<\lambda(x_0)/2,\;n\geq0.                                                      \vspace{-0.2cm}\eqno(3.2)$$

Suppose that system (2.1) is not strongly sensitive at $x_0$. Then,
for any $0<\eta<\delta_{1}$, there exists $y_0\in I$ with $0<|y_0-x_0|<\eta$ such that
\vspace{-0.2cm}$$|w_n|=|y_n-x_n|\leq\eta,\;n>0,                                                                                \vspace{-0.2cm}\eqno(3.3)$$
where $\{x_n\}_{n=0}^{\infty}$ and $\{y_n\}_{n=0}^{\infty}$ are the orbits of system (2.1) starting from  $x_0$ and $y_0$, respectively,
and
\vspace{-0.2cm}$$w_n:=y_n-x_n,\;n\geq0.                                                                                        \vspace{-0.2cm}\eqno(3.4)$$
It is clear that one can inductively get that $w_n\neq0$ for each $n\geq0$.
So, there exists $c_n\in(x_n, y_n)$ (or $(y_n, x_n)$) such that
\vspace{-0.2cm}$$w_{n+1}=f_n(x_n+w_n)-f_n(x_n)=f_n'(c_n)w_n,\;n\geq0.                                                          \vspace{-0.2cm}\eqno(3.5)$$
Then $|c_n-x_n|\leq |w_n|\leq\eta<\delta_{1}<\delta_{0}$,
which implies that $c_{n}\in J$ for each $n\geq0$. It follows from (3.2) that for each $n\geq0$,
\vspace{-0.2cm}$$\big|\frac{1}{n+1}\sum_{k=0}^{n}\ln|f_k'(c_k)|-\frac{1}{n+1}\sum_{k=0}^{n}\ln|f_k'(x_k)|\big|
\leq\frac{1}{n+1}\sum_{k=0}^{n}\big|\ln|f_k'(c_k)|-\ln|f_k'(x_k)|\big| <\lambda(x_0)/2.\vspace{-0.2cm}$$
Then
\vspace{-0.2cm}$$\frac{1}{n+1}\sum_{k=0}^{n}\ln|f_k'(c_k)|>\frac{1}{n+1}\sum_{k=0}^{n}\ln|f_k'(x_k)|
-\lambda(x_0)/2,\;n\geq0.\vspace{-0.2cm}$$
So
\vspace{-0.2cm}$$\limsup_{n\to\infty}\frac{1}{n+1}\sum_{k=0}^{n}\ln|f_k'(c_k)|\geq\lambda(x_0)-\lambda(x_0)/2
=\lambda(x_0)/2,\vspace{-0.2cm}$$
which yields that there exists an increasing sequence $\{n_j\}_{j=0}^{\infty}$ such that
\vspace{-0.2cm}$$\lim_{j\to\infty}\frac{1}{n_j+1}\sum_{k=0}^{n_j}\ln|f_k'(c_k)|
\geq\lambda(x_0)/2.\vspace{-0.2cm}$$
Then there exists a positive integer $k_{0}$ such that for each $j\geq k_{0}$, one has that
\vspace{-0.3cm}$$\frac{1}{n_j+1}\sum_{k=0}^{n_j}\ln|f_k'(c_k)|\geq\lambda(x_0)/4.        \vspace{-0.2cm}\eqno(3.6)$$
It follows from (3.5) and (3.6) that for each $j\geq k_{0}$,
\vspace{-0.2cm}$$|w_{n_{j}+1}|=|f_{n_j}'(c_{n_j})|\cdots|f_{0}'(c_{0})||w_{0}|
\geq\exp[(n_{j}+1)\lambda(x_0)/4]|w_{0}|,                                              \vspace{-0.2cm}\eqno(3.7)$$
then $|w_{n_{j}+1}|\to\infty$ as $j\to\infty$ since $\lambda(x_0)>0$. This is a contradiction with (3.3).
Therefore, system (2.1) is strongly sensitive at $x_0$.
This completes the proof.\medskip

\noindent{\bf Remark 3.1.} The result of Theorem 3.1 extends that of [16, Theorem 1] for autonomous discrete
systems to non-autonomous discrete systems.
Moreover, the condition about the existence of the limit of $\sum_{k=0}^{n-1}\ln|f_{k}'(x_k)|/n$
in [16, Theorem 1] is removed in Theorem 3.1.
\medskip

The following result shows that the positive Lyapunov exponent in a totally invariant set implies sensitivity in this set
under certain conditions.\medskip

\noindent{\bf Theorem 3.2.} {\it Let $I$ be a non-degenerate interval, $f_n: I\to I$ be a $C^{1}$ map for each $n\geq0$,
and $\Lambda$ be a totally invariant subinterval in $I$; that is, $f_n(\Lambda)\subset\Lambda,\;n\geq0$. Assume that
$\{f_{n}'\}_{n=0}^{\infty}$ is equi-continuous in $\Lambda$, and there exists $M>0$ such that for any $x\in\Lambda$,
$|f_{n}'(x)|\geq M,\;n\geq0$. If $\lambda:=\inf_{x\in\Lambda}\lambda(x)>0$,
then system {\rm(2.1)} is strongly sensitive in $\Lambda$.}\medskip

\noindent{\bf Proof.} With a similar argument to that used in the proof of Theorem 3.1, one can show this theorem. The proof is complete.\bigskip

\noindent{\bf 4. Negative Lyapunov exponent implies exponential asymptotical stability}\medskip

In this section, we shall investigate some relationship between negative Lyapunov exponents and stability
for a class of non-autonomous discrete systems.\medskip

\noindent{\bf Lemma 4.1} [16, Lemma 4]. {\it Let $z_n$, $\mu_n (n\geq 0)$, and $B$ are nonnegative numbers. If
$z_n\leq B+\sum_{k=1}^{n}\mu_{k}z_{k-1}$ for each $n\geq0$, then
$z_n\leq B\exp{(\sum_{k=1}^{n}\mu_{k})}$ for each $n\geq0$.}\medskip

Denote\vspace{-0.2cm}
$$\lambda_{0}(x_0):=\liminf_{n\to\infty}\frac{1}{n}\sum_{k=0}^{n-1}\ln|f_{k}'(x_k)|. \vspace{0.5mm}    \eqno(4.1) $$

\noindent{\bf Theorem  4.1.} {\it Let $I$ be a non-degenerate, closed, and bounded interval, $f_n: I\to I$ be a $C^{2}$ map
for each $n\geq0$, and $x_0\in I$. Assume that $\{f_{n}''\}_{n=0}^{\infty}$ is uniformly bounded in $I$.
If $-\infty<\lambda(x_0)<0$ and $2\lambda(x_0)<\lambda_{0}(x_0)$,
then system {\rm(2.1)} is exponentially asymptotically stable at $x_0$.}\medskip

\noindent{\bf Proof.} For each $n\geq0$ and any $y_0\in I$ with $y_0\neq x_0$, set
\vspace{-0.2cm}$$w_{n+1}=a_{n}w_n+g_n(w_n),                                                                                        \vspace{-0.2cm}\eqno(4.2)$$
where $w_n$ is specified in (3.4) and
\vspace{-0.2cm}$$a_n:=f'_n(x_n),\;\;g_n(w_n):=f_n(x_n+w_n)-f_n(x_n)-f_n'(x_n)w_n.                                                  \vspace{-0.2cm}\eqno(4.3)$$
By (4.2) one can inductively get that
\vspace{-0.4cm}$$w_n=a_{n-1}\cdots a_0w_0+\sum_{k=1}^{n}\big(a_{n-1}\cdots a_kg_{k-1}(w_{k-1})\big),\;n\geq1.                              \vspace{-0.2cm}\eqno(4.4)$$
Since $\{f_{n}''\}_{n=0}^{\infty}$ is uniformly bounded in $I$, there exists $M>0$ such that $|f_n''(x)|\leq M$ for every $x\in I$ and each $n\geq0$.
This, together with the second relation in (4.3), yields that
\vspace{-0.2cm}$$|g_n(w_n)|\leq\frac{1}{2}M|w_n|^{2},\;n\geq0.                                                                      \eqno(4.5)\vspace{-0.2cm}$$
It follows from (2.4) and (4.1) that for any fixed $0<\epsilon_{0}<[\lambda_{0}(x_0)-2\lambda(x_0)]/3$,
there exists a positive integer $N_{0}$ such that
\vspace{-0.2cm}$$\lambda_{0}(x_0)-\epsilon_{0}<\frac{1}{n}\sum_{i=0}^{n-1}\ln|a_i|<\lambda(x_0)+\epsilon_{0},\; n>N_{0}.\vspace{-0.2cm}$$
Then, for any $n>k>N_{0}$, one has that
\vspace{-0.2cm}$$|a_{n-1}\cdots a_k|=\frac{|a_{n-1}\cdots a_0|}{|a_{k-1}\cdots a_0|}<\exp{[\lambda(n-k)+l_k]},\vspace{-0.2cm}$$
where $\lambda:=\lambda(x_0)+\epsilon_{0}<0$, and $l_{k}:=k[\lambda(x_0)-\lambda_{0}(x_0)+2\epsilon_{0}]$.
Note that $|a_i|>0$ for each $i\geq0$ since $\lambda(x_0)>-\infty$. Because $f_n'$
is continuous in the closed and bounded interval $I$ for each $n\geq0$,
it can be easily verified that there exists a constant $C_{0}\geq1$ depending on
$\epsilon_{0}$ and $N_{0}$ such that for each $n\geq k\geq0$,
\vspace{-0.2cm}$$|a_{n-1}\cdots a_k|\leq C_{0}\exp{[\lambda(n-k)+l_k]}.                                                             \vspace{-0.2cm}\eqno(4.6)$$
Then, by (4.4)--(4.6) one has that
\vspace{-0.4cm}$$|w_n|\leq C_{0}\exp{[\lambda n]}|w_0|+\frac{1}{2}MC_{0}\sum_{k=1}^{n}\exp{[\lambda(n-k)+l_{k}]}|w_{k-1}|^{2},\;n\geq0. \vspace{-0.2cm}\eqno(4.7)$$

Next, for any $\eta>0$, we shall inductively show that
\vspace{-0.2cm}$$|w_n|=|y_n-x_n|\leq\eta\exp{[\lambda n]},\;n\geq0,                      \vspace{-0.2cm}\eqno(4.8)$$
for any $y_0\in I$ with $|y_0-x_0|<\delta$, where
$$\delta:= C_{0}^{-1}\exp{[-D_{\eta}]}\eta,\;\;
D_{\eta}:=\frac{1}{2}MC_{0}\eta\exp{[-2\lambda]}\exp{[\widetilde{\lambda}]}
\left(1-\exp{[\widetilde{\lambda}]}\right)^{-1}$$
with $\widetilde{\lambda}:=2\lambda(x_0)-\lambda_{0}(x_0)+3\epsilon_0<0$.

Evidently, (4.8) holds for $n=0$. Suppose that (4.8) holds for each $0\leq n\leq T-1$,
where $T$ is a positive integer. We show that there exists $\delta>0$, which only depends on $\eta$ but not on $T$,
such that (4.8) holds for $n=T$. By the assumption, $|w_{k-1}|<\eta\exp{[\lambda (k-1)]}$ for each
$1\leq k\leq n\leq T$. This, together with (4.7), implies that
\vspace{-0.2cm}$$|w_n|\leq C_{0}\exp{[\lambda n]}|w_0|+\frac{1}{2}MC_{0}\eta\sum_{k=1}^{n}\exp{[\lambda(n-1)+l_{k}]}|w_{k-1}|,\;1\leq n\leq T.
                                                                                   \vspace{-0.2cm} \eqno(4.9)$$
Let $z_{n}=|w_n|\exp{[-\lambda n]}$. It follows from (4.9) that
\vspace{-0.2cm}$$z_n\leq C_{0}|w_0|+\frac{1}{2}MC_{0}\eta\sum_{k=1}^{n}
\exp{[\lambda(k-2)+l_{k}]}z_{k-1},\;1\leq n\leq T.\vspace{-0.2cm}$$
Thus by Lemma 4.1 we get that
\vspace{-0.2cm}$$z_n\leq C_{0}|w_0|\exp{\left\{\frac{1}{2}MC_{0}\eta\exp{[-2\lambda]}
\sum_{k=1}^{n}\exp{[\lambda k+l_{k}]}\right\}},\;1\leq n\leq T.                          \vspace{-0.2cm}\eqno(4.10)$$
Since $\lambda k+l_{k}=k[2\lambda(x_0)-\lambda_{0}(x_0)+3\epsilon_0]=k\widetilde{\lambda}$,
$\sum_{k=1}^{n}\exp{[\lambda k+l_{k}]}=\sum_{k=1}^{n}\exp{[k\widetilde{\lambda}]}
<\exp{[\widetilde{\lambda}]}(1-\exp{[\widetilde{\lambda}]})^{-1}$. Thus, by (4.10) one has that
\vspace{-0.2cm}$$|w_n|=z_n\exp{[\lambda n]}\leq C_{0}|w_0|\exp{[D_{\eta}]}\exp{[\lambda n]},\; 1\leq n\leq T.\vspace{-0.1cm}$$
In particular,
$|w_T|\leq C_{0}|w_0|\exp{[D_{\eta}]}\exp{[\lambda T]}<\eta\exp{[\lambda T]}$.
Hence, (4.8) also holds for $n=T$.
Therefore, system (2.1) is exponentially asymptotically stable at $x_0$. The proof is complete.\medskip

If the limit of $(\sum_{k=0}^{n-1}\ln|f_{k}'(x_k)|)/n$ exists, then $\lambda(x_0)=\lambda_{0}(x_0)$ and the following
result can be directly derived from Theorem 4.1.\medskip

\noindent{\bf Corollary 4.1.} {\it Let $I$ be a non-degenerate, closed, and bounded interval, $f_n: I\to I$ be a $C^{2}$ map
for each $n\geq0$, and $x_0\in I$. Assume that $\{f_{n}''\}_{n=0}^{\infty}$ is uniformly bounded in $I$.
If $-\infty<\lambda(x_0)=\lambda_{0}(x_0)<0$, then system {\rm(2.1)} is exponentially asymptotically stable
at $x_0$.}\medskip

\noindent{\bf Remark 4.1.} The result of  Corollary 4.1 extends that of [16, Theorem 4]
for autonomous systems to non-autonomous systems.\bigskip

\noindent{\bf 5. An example}\medskip

Consider the following non-autonomous logistic system:
\vspace{-0.2cm}$$x_{n+1}=r_{n}x_{n}(1-x_{n}),\;\; n\geq0,                                                                      \vspace{-0.2cm}\eqno(5.1)$$
governed by the maps $f_{n}(x)=r_{n}x(1-x)$, $x\in I:=[0,1]$. Let $0< r_n\leq4,\;n\geq0$. Then $f_n(I)\subset I$ for each $n\geq0$.

It is evident that for each $n\geq0$, $f_{n}$ is $C^{2}$ in $I$, $f_{n}'(x)=r_{n}(1-2x)$, and $f_n''(x)=-2r_n$.
Clearly, $0$ is a fixed point of system (5.1) and $f_{n}'(0)=r_n,\;n\geq0$,
which results in
\vspace{-0.2cm}$$\lambda(0)=\limsup_{n\to\infty}\frac{1}{n}\sum_{k=0}^{n-1}\ln|f_{k}'(0)|=\limsup_{n\to\infty}\frac{1}{n}\sum_{k=0}^{n-1}\ln r_{k}.                                                                                       \vspace{-0.2cm}\eqno(5.2)$$

The limit of $(\sum_{k=0}^{n-1}\ln r_{k})/n$ may not exist. However, its upper limit exists. Hence, the upper limit used
in Definition 2.6 is reasonable in this case.

In the case that $L\leq r_{n}\leq4,\;n\geq0$, where $L>1$ is a constant, one has that
$|f_n'(0)|\geq L>1$ and $|f_n''(x)|=2r_n\leq8$ for all $n\geq0$.
So $\{f_n'\}_{n=0}^{\infty}$ is equi-continuous in $I$. It follows from (5.2) that $\lambda(0)\geq\ln L>0$.
Hence, all the assumptions in Theorem 3.1 hold for system (5.1) with $x_0=0$.
Therefore, system (5.1) is strongly sensitive at $0$ in this case.

In the case that $0<a\leq r_{n}\leq b<1,\;n\geq0$, where $a$ and $b$ are positive constants with $b^{2}<a$,
by (5.2) we get that
\vspace{-0.3cm}$$\lambda(0)\leq\ln b<0, \;\lambda_{0}(0)\geq \ln a>-\infty.                                                    \vspace{-0.2cm}\eqno(5.3)$$
Since $b^{2}<a$, it follows from (5.3) that
\vspace{-0.2cm}$$2\lambda(0)-\lambda_{0}(0)\leq2\ln b-\ln a=\ln(b^{2}/a)<0.                                                    \vspace{-0.2cm}$$
Moreover, $|f_n''(x)|=2r_n<2$ for each $n\geq0$, which yields that $\{f_n''\}_{n=0}^{\infty}$ is uniformly bounded in $I$.
Hence, all the assumptions in Theorem 4.1 hold for system (5.1) with $x_0=0$. Therefore, system (5.1)
is exponentially asymptotically stable at $0$.
\bigskip

\noindent{\bf Acknowledgment}\medskip

This research was supported by the NNSF of China (grant 11571202).
\end{document}